\documentclass[12pt,a4paper]{article}
\usepackage{mathrsfs}
\usepackage{amsfonts}
\usepackage{}
\usepackage[top=2.5cm,bottom=2.5cm,left=2.5cm,right=2.5cm]{geometry}
\def\Dj{\hbox{D\kern-.73em\raise.30ex\hbox{-}
\raise-.30ex\hbox{}}}
\def\dj{\hbox{d\kern-.33em\raise.80ex\hbox{-}
\raise-.80ex\hbox{\kern-.40em}}}
\usepackage{epsfig}
\usepackage{amsmath,amsthm,amsfonts,amssymb,amscd,cite,graphicx,array}

\newtheorem{thm}{Theorem}

\newtheorem{lem}{Lemma}
\newtheorem{ob}{Observation}

\newtheorem{cor}{Corollary}
\newtheorem{pro}{Proposition}

\def\qed{\hfill \nopagebreak\rule{5pt}{8pt}}
\def\pf{\noindent {\it Proof.} }

\makeatletter \@addtoreset{equation}{section} \makeatother

\begin{document}

\begin{center} {\Large \bf The $k$-proper index of graphs\footnote{Supported by NSFC No. 11371205, 11531011.}}
\end{center}

\begin{center}
{{ \small Lin Chen, Xueliang Li, Jinfeng Liu} \\[3mm]
{\small Center for Combinatorics and LPMC}\\
{\small Nankai University, Tianjin 300071, P.R. China}\\[3mm]
{\small E-mail: chenlin1120120012@126.com, lxl@nankai.edu.cn, ljinfeng709@163.com}}
\end{center}

\begin{center}
\begin{minipage}{130mm}
\vskip 0.2cm
\begin{center}
{\bf Abstract}
\end{center}
{\small A tree $T$ in an edge-colored graph is a \emph{proper tree} if any two adjacent edges of $T$ are colored with different colors. Let $G$ be a graph of order $n$ and $k$ be a fixed integer with $2\leq k\leq n$. For a vertex set $S\subseteq V(G)$, a tree containing the vertices of $S$ in $G$ is called an \emph{$S$-tree}. An edge-coloring of $G$ is called a \emph{$k$-proper coloring} if for every set $S$ of $k$ vertices in $G$, there exists a proper $S$-tree in $G$. The \emph{$k$-proper index} of a nontrivial connected graph $G$, denoted by $px_k(G)$, is the smallest number of colors needed in a $k$-proper coloring of $G$. In this paper, some simple observations about $px_k(G)$ for a nontrivial connected graph $G$ are stated. Meanwhile, the $k$-proper indices of some special graphs are determined, and for every pair of positive integers $a$, $b$ with $2\leq a\leq b$, a connected graph $G$ with $px_k(G)=a$ and $rx_k(G)=b$ is constructed for each integer $k$ with $3\leq k\leq n$. Also, the graphs with $k$-proper index $n-1$ and $n-2$ are respectively characterized.\\
[2mm] {\bf Keywords:} coloring of graphs, $k$-proper index, characterization of graphs \\
[2mm] {\bf AMS Subject Classification (2010):} 05C05, 05C15, 05C40. }
\end{minipage}
\end{center}

\section{Introduction}
In this paper, all graphs under our consideration are finite, undirected, connected and simple. For more notation and terminology that will be used in the sequel, we refer to \cite{graphBondy2008}, unless otherwise stated.

In 2008, Chartrand et al. \cite{CJMZ} first introduced the concept of rainbow connection. Let $G$ be a nontrivial connected graph on which an edge-coloring $c: E(G)\rightarrow \{1,2,\ldots,k\}\\
(k\in \mathbb{N})$ is defined, where adjacent edges may be colored with the same color. For any two vertices $u$ and $v$ of $G$, a path in $G$ connecting $u$ and $v$ is abbreviated as a $uv$-path. A $uv$-path $P$ is a \emph{rainbow $uv$-path} if no two edges of $P$ are colored with the same color. The graph $G$ is \emph{rainbow connected} (with respect to $c$) if $G$ contains a rainbow $uv$-path for every two vertices $u$ and $v$, and the coloring $c$ is called a \emph{rainbow coloring} of $G$. If $k$ colors are used, then $c$ is a \emph{rainbow $k$-coloring}. The minimum $k$ for which there exists a rainbow $k$-coloring of the edges of $G$ is the \emph{rainbow connection number} of $G$, denoted by $rc(G)$. The topic of rainbow connection is fairly interesting and numerous relevant papers have been written. For more details see a survey \cite{LSS} and a book \cite{LS}.

Subsequently, a series of generalizations of rainbow connection number were proposed. The $k$-rainbow index is one of them. An edge-colored tree $T$ is a \emph{rainbow tree} if no two edges of $T$ are assigned the same color. Let $G$ be a nontrivial connected graph of order $n$ and let $k$ be an integer with $2\leq k\leq n$. A \emph{$k$-rainbow coloring} of $G$ is an edge coloring of $G$ having the property that for every set $S$ of $k$ vertices of $G$, there exists a rainbow tree $T$ in $G$ such that $S\subseteq V(T)$. The minimum number of colors needed in a $k$-rainbow coloring of $G$ is the \emph{$k$-rainbow index} of $G$. These concepts were introduced by Chartrand et al. in \cite{COZ}, and were further studied in \cite{CLZ,CLZ2,CLYZ,LSYZ,LSYZ2,LiuH}.

In addition, a natural extension of the rainbow connection number is the proper connection number, which was introduced by Borozan et al. in \cite{BFGMMMT}. A path in an edge-colored graph is said to be \emph{properly edge-colored (or proper)}, if every two adjacent edges differ in color. An edge-colored graph $G$ is \emph{$k$-proper connected} if any two vertices are connected by $k$ internally pairwise vertex-disjoint proper paths. The \emph{$k$-proper connection number} of a $k$-connected graph $G$, denoted by $pc_k(G)$, is defined as the smallest number of colors that are needed in order to make $G$ $k$-proper connected. In particular, when $k=1$, the $1$-proper connection number is abbreviated as proper connection number and written as $pc(G)$. For more results, we refer to \cite{ALLZ,FGM,GLQ,HLW1,HLW2,LLZ,LWY}.

Inspired by the $k$-rainbow index and the proper connection number, a natural idea is to introduce the concept of $k$-proper index. A tree $T$ in an edge-colored graph is a \emph{proper tree} if any two adjacent edges of $T$ are colored with different colors. Let $G$ be a graph of order $n$ and $k$ be a fixed integer with $2\leq k\leq n$. For a vertex set $S\subseteq V(G)$, a tree containing the vertices of $S$ in $G$ is called an \emph{$S$-tree}. An edge-coloring of $G$ is called a \emph{$k$-proper coloring} if for every set $S$ of $k$ vertices in $G$, there exists a proper $S$-tree in $G$. The \emph{$k$-proper index} of a nontrivial connected graph $G$, denoted by $px_k(G)$, is the smallest number of colors needed in a $k$-proper coloring of $G$. By definition, $px_2(G)$ is precisely the proper connection number $pc(G)$ for any nontrivial graph $G$. As a variety of nice results about $pc(G)=px_2(G)$ have been obtained, we in this paper only study $px_k(G)$ for $3\leq k\leq n$.

The paper is organized as follows: In Section \ref{Preliminaries}, some simple observations about $px_k(G)$ for a nontrivial graph $G$ are stated. Meanwhile, certain necessary lemmas are also listed. In Section \ref{Specialgraphs}, the $k$-proper indices of some special graphs are determined. And for every pair of positive integers $a$, $b$ with $2\leq a\leq b$, a connected graph $G$ with $px_k(G)=a$ and $rx_k(G)=b$ is constructed for each integer $k$ with $3\leq k\leq n$. In Section \ref{Characterize}, the graphs with $k$-proper index $n-1$ and $n-2$ are characterized, respectively.

\section{Preliminaries}\label{Preliminaries}
We in this section state some observations about $px_k(G)$ for a nontrivial graph $G$. Also, certain necessary lemmas are listed.

For a graph $G$ with order $n\geq 3$, it follows from the definition that
$$pc(G)=px_2(G)\leq px_3(G)\leq px_4(G)\leq \cdots\leq px_n(G). \leqno(*)$$
This simple property will be used frequently later.

Since any $k$-proper coloring of a spanning subgraph must be a $k$-proper coloring of its supergraph. Then there exists a fundamental proposition about spanning subgraphs.

\begin{pro}\label{subgraph}
If $G$ is a nontrivial connected graph of order $n\geq 3$ and $H$ is a connected spanning subgraph of $G$, then $px_k(G)\leq px_k(H)$ for any $k$ with $3\leq k\leq n$. In particular, $px_k(G)\leq px_k(T)$ for every spanning tree $T$ of $G$.
\end{pro}

It has been seen in \cite{COZ} that $rx_k(G)\leq n-1$ for any graph $G$ with order $n\geq 3$ and any integer $k$ with $3\leq k\leq n$. Since a rainbow tree must be a proper tree, then obviously $px_k(G)\leq rx_k(G)\leq n-1$. Moreover, this simple upper bound is sharp, the graphs with $px_k(G)=n-1$ will be characterized in Section \ref{Characterize}.

For any nontrivial graph $G$, $\chi'(G)$ denotes the edge-chromatic number of $G$. It is well-known that either $\chi'(G)=\Delta(G)$ or $\chi'(G)=\Delta(G)+1$ by Vizing's Theorem, where $\Delta(G)$, or simply $\Delta$, is the maximum degree of $G$. Accordingly, a natural upper bound of $px_k(G)$ with respect to these parameters follows.

\begin{pro}\label{upperbound}
Let $G$ be a graph with order $n\geq 3$, maximum degree $\Delta(G)$ and edge-chromatic number $\chi'(G)$. Then for each integer $k$ with $3\leq k\leq n$, we have
$$px_k(G)\leq \chi'(G)\leq \Delta(G)+1.$$
\end{pro}

In addition, there is a classical result about the edge-chromatic number of a graph, which will be useful in the next section.

\begin{lem}[\cite{graphBondy2008}]\label{edgechromaticnumber}
If $G$ is bipartite, then $\chi'(G)=\Delta(G)$.
\end{lem}

For arbitrary $k \ (k\geq 3)$ vertices of a nontrivial graph $G$, any tree $T$ containing these vertices must contain internal vertices. While for any proper tree $T$, there must be $d(u)$ distinct colors assigned to the incident edges of each vertex $u$ in $T$, where $d(u)$ denotes the degree of $u$ in $T$. Hence, the incident edges of any internal vertex must be assigned with at least two distinct colors to make $T$ proper. Then the following trivial lower bound is immediate.

\begin{pro}\label{lowerbound}
For arbitrary graph $G$ with order $n\geq 3$, we have
$$px_k(G)\geq 2$$
for any integer $k$ with $3\leq k\leq n$.
\end{pro}

{\bf Remark:} The above lower bound of $px_k(G)$ is sharp since there exist many graphs satisfying $px_k(G)=2$, as shown in Section \ref{Specialgraphs}. Further, we believe that it will be interesting to characterize all graphs with $k$-proper index 2 for specific values of $k$.

In any graph $G$, a path (resp. cycle) that contains every vertex of $G$ is called a \emph{Hamilton path} (resp. \emph{Hamilton cycle}) of $G$. A graph is \emph{traceable} if it contains a Hamilton path, and a graph is \emph{hamiltonian} if it contains a Hamilton cycle. The following is an immediate consequence of these definitions as well as Proposition \ref{lowerbound}.

\begin{pro}\label{traceable}
If $G$ is a traceable graph with order $n\geq 3$, then $px_k(G)=2$ for each integer $k$ with $3\leq k\leq n$.
\end{pro}

As mentioned before, characterizing all graphs with $k$-proper index 2 for specific values of $k$ would be interesting.
While for the cases of $k=n$ and $k=n-1$, there are two basic results that can be mentioned.

\begin{ob}
If a graph $G$ of order $n$ satisfies $px_n(G)=2$, namely, $px_k(G)=2$ for each $k$ with $3\leq k\leq n$. Then $G$ is a traceable graph.
\end{ob}

\begin{ob}
If a graph $G$ of order $n$ satisfies $px_k(G)=2$ for each $k$ with $3\leq k\leq n-1$. Then $px_n(G)=2$ if and only if $G$ is traceable. Otherwise, $px_n(G)=3$.
\end{ob}

It is well known that if $G$ is a simple graph with order $n\geq 3$ and minimum degree $\delta\geq \frac{n}{2}$, then $G$ is hamiltonian. Whereupon a direct corollary follows.

\begin{cor}
If $G$ is a simple graph with order $n\geq 3$ and minimum degree $\delta\geq \frac{n}{2}$, then $px_k(G)=2$ for each integer $k$ with $3\leq k\leq n$.
\end{cor}

In \cite{COZ}, Chartrand et al. derived the $k$-rainbow index of a nontrivial tree, which will be helpful in the next section.

\begin{lem}[\cite{COZ}]\label{rtrees}
Let $T$ be a tree of order $n\geq 3$. For each integer $k$ with $3\leq k\leq n$,
$$rx_k(T)=n-1.$$
\end{lem}

In \cite{BFGMMMT}, Borozan et al. established the proper connection number of trees.

\begin{lem}[\cite{BFGMMMT}]\label{pc(T)}
If $G$ is a tree then $pc(G)=\Delta(G)$.
\end{lem}

At the end of this section, we recall several notations required in the subsequent sections.

Let $E'\subseteq E(G)$ be a set of edges of a graph $G$, then $G[E']$ is the subgraph of $G$ induced by $E'$. If $e$ is an edge of $G$, then $G-e$ denotes the graph obtained from $G$ by only deleting the edge $e$. If $G$ is not complete, denote by $G+e$ the graph obtained from $G$ by the addition of $e$, where $e$ is an edge connecting two nonadjacent vertices of $G$.

\section{The $k$-proper indices of special graphs}\label{Specialgraphs}
In this section, we determine the $k$-proper indices of complete graphs, cycles, wheels, trees and unicyclic graphs. Moreover, the independence of $px_k(G)$ and $rx_k(G)$ is given by a brief theorem.

It has been seen that if $G$ is a traceable graph, then $px_k(G)=2$. Obviously, the complete graphs, cycles and wheels are all traceable, thus the $k$-proper indices of these graphs are direct consequences of Proposition \ref{traceable}.

\begin{thm}\label{completegraphs}
Let $K_n$, $C_n$ and $W_n$ be a complete graph, a cycle and a wheel with $n \ (n\geq 3)$ vertices, respectively. Then for any integer $k$ with $3\leq k\leq n$, we have
$$px_k(K_n)=px_k(C_n)=px_k(W_n)=2.$$
\end{thm}

Now we determine the $k$-proper index for a nontrivial tree.

\begin{thm}\label{trees}
If $T$ is a tree of order $n\geq 3$, then for each integer $k$ with $3\leq k\leq n$,
$$px_k(T)=\Delta(T).$$
\end{thm}
\pf Firstly, since $T$ is bipartite, then $px_k(T)\leq \chi'(T)=\Delta(T)$ for $3\leq k\leq n$ by Proposition \ref{upperbound} and Lemma \ref{edgechromaticnumber}. On the other hand, according to Ineq. $(*)$ and Lemma \ref{pc(T)}, $px_k(T)\geq pc(T)=\Delta(T)$ holds naturally for $3\leq k\leq n$. Therefore, we arrive at $px_k(T)=\Delta(T)$ for any integer $k$ with $3\leq k\leq n$. \qed

Combine with Proposition \ref{subgraph} and Theorem \ref{trees}, one can check that the following proposition holds.

\begin{pro}\label{spanning}
For any graph $G$ with order $n\geq 3$ and any integer $k$ with $3\leq k\leq n$, we have
$$px_k(G)\leq \min\{\Delta(T): T~\text{is a spanning tree of}~ G\}.$$
\end{pro}

Since $\Delta(T)\leq \Delta(G)$ for any spanning tree $T$ of $G$. Then the upper bound in Proposition \ref{upperbound} can be replaced by $\Delta(G)$.

\begin{pro}\label{upperbound2}
Let $G$ be a graph with order $n\geq 3$ and maximum degree $\Delta(G)$, then
$$px_k(G)\leq \Delta(G)$$
for each integer $k$ with $3\leq k\leq n$.
\end{pro}

{\bf Remark:} The above upper bound of $px_k(G)$ is sharp since the equality holds apparently for arbitrary nontrivial tree.

In order to get the $k$-proper index of a unicyclic graph, an assistant lemma is presented.

\begin{lem}\label{bridge}
Let $G$ be a graph of order $n\geq 3$ containing bridges and $v$ be any vertex of $G$. Denote by $b(v)$ the number of bridges incident with $v$. Set $b=\max\{b(v): v\in V(G)\}$. Then for each integer $k$ with $3\leq k\leq n$, we have $px_k(G)\geq b$.
\end{lem}
\pf Since for $3\leq k\leq n$, it has been seen from Ineq. $(*)$ that $px_k(G)\geq px_3(G)$. Then we should only prove the case when $k=3$. Since $px_3(G)\geq 2$ by Proposition \ref{lowerbound}, the result is trivial when $b=1$ or $2$. Thus we may assume that $b\geq 3$. Suppose that $b(u)=b=\max\{b(v): v\in V(G)\}$ for some vertex $u$. Let $F=\{uw_1,uw_2,\ldots,uw_b\}$ be the set of bridges incident with $u$. Set $A=\{u,w_1,w_2,\ldots,w_b\}$. For any $3$-set $S=\{w_i,w_j,u\}\subseteq A$, where $i,j\in\{1,2,\ldots,b\}$ and $i\neq j$, every $S$-tree $T$ must contain the edges $uw_i$ and $uw_j$. Hence, the edges $uw_i$ and $uw_j$ receive distinct colors to make $T$ proper. Which implies that the edges $uw_1,uw_2,\ldots,uw_b$ need $b$ distinct colors in any $3$-proper coloring of $G$. Therefore, $px_3(G)\geq b$. This completes the proof. \qed

With the aid of Lemma \ref{bridge}, now we are able to deal with the $k$-proper index for a unicyclic graph.

\begin{thm}\label{unicyclic}
Let $G$ be a unicyclic graph of order $n\geq 3$, and maximum degree $\Delta(G)$. Then, for each integer $k$ with $3\leq k\leq n$,
$$px_k(G)=\Delta(G)-1$$
when $G$ contains at most two vertices having maximum degree such that the vertices with maximum degree are all in the unique cycle of
$G$ and these two vertices (if both exist) are adjacent;\\
Otherwise,
$$px_k(G)=\Delta(G).$$
\end{thm}
\pf Note that when $G=C_n$, it follows from Theorem \ref{completegraphs} that $px_k(G)=px_k(C_n)=2=\Delta(G)$ for $3\leq k\leq n$. Thus in the following we assume that $G$ is not a cycle. And assume the vertices in the unique cycle of $G$ are $u_1,u_2,\ldots,u_g$. Besides, keep in mind that $px_k(G)\leq \Delta(G)$ for $3\leq k\leq n$, which will be used later. As before, denote by $b(v)$ the number of bridges incident with the vertex $v$. The discussion is divided into three cases.

{\bf Case 1.} At first, assume that $G$ contains a vertex, say $u$, satisfying\\
(1) the degree of $u$ is $d(u)=\Delta(G)$.\\
(2) $u$ is not in the cycle of $G$.

Then evidently the incident edges of $u$ are all bridges, i.e., $b(u)=d(u)=\Delta(G)$. According to Lemma \ref{bridge}, we directly have $px_k(G)\geq b(u)=\Delta(G)$ for $3\leq k\leq n$. Meanwhile, Proposition \ref{upperbound2} guarantees $px_k(G)\leq \Delta(G)$ for $3\leq k\leq n$. Accordingly, we get $px_k(G)=\Delta(G)$ for each integer $k$ with $3\leq k\leq n$ in this case.

By Case 1, if such a vertex $u$ exists in $G$, then we always have $px_k(G)=\Delta(G)$ for each integer $k$ with $3\leq k\leq n$. To avoid redundant presentation, we in the following suppose that $G$ contains no such vertices.

{\bf Case 2.} Now assume $G$ simultaneously satisfies\\
(3) $G$ contains at most two vertices having maximum degree;\\
(4) the vertices with maximum degree are all in the unique cycle of $G$;\\
(5) these two vertices (if both exist) are adjacent in $G$.

Then without loss of generality, suppose that $d(u_1)=\Delta(G)$, $d(u_2)\leq \Delta(G)$ and $d(u)<\Delta(G)$ for any other vertex $u$. Moreover, suppose that the neighbors of $u_1$ are $v_1$, $v_2$, $\ldots$, $v_{\Delta(G)-2}$, $v_{\Delta(G)-1}=u_2$ and $v_{\Delta(G)}=u_g$. Thereupon, in any $3$-proper coloring $c$ of $G$, based on the proof of Lemma \ref{bridge}, the edges $u_1v_i$ with $i\in \{1,2,\ldots,\Delta(G)-2\}$ are assigned with $\Delta(G)-2$ distinct colors since they are all bridges incident with $u_1$. Without loss of generality, suppose that $c(u_1v_1)=1$, $c(u_1v_2)=2$, $\ldots$, $c(u_1v_{\Delta(G)-2})=\Delta(G)-2$. Further, we claim that at least a new color is used by the edges $u_1u_2$ and $u_1u_g$. For otherwise, suppose that $c(u_1u_2)=i$ and $c(u_1u_g)=j$ with $i,j\in \{1,2,\ldots,\Delta(G)-2\}$. If $i=j$, then there exists no proper tree containing the vertices $u_1$, $u_2$ and $v_i$, a contradiction. If $i\neq j$, then there exists no proper tree containing the vertices $v_i$, $v_j$ and $u_2$, again a contradiction. Therefore, at least $\Delta(G)-2+1=\Delta(G)-1$ different colors are used by $c$. It follows that $px_3(G)\geq \Delta(G)-1$. Thus, Ineq.$(*)$ deduces that $px_k(G)\geq px_3(G)\geq \Delta(G)-1$ for each integer $k$ with $3\leq k\leq n$. On the other hand, obviously $G-u_1u_2$ is a spanning tree of $G$ with maximum degree $\Delta(G)-1$. By Theorem \ref{trees}, we know that $px_k(G-u_1u_2)=\Delta(G-u_1u_2)=\Delta(G)-1$ for $3\leq k\leq n$. Consequently, $px_k(G)\leq px_k(G-u_1u_2)=\Delta(G)-1$ based on Proposition \ref{subgraph}. To sum up, we obtain $px_k(G)=\Delta(G)-1$ for each integer $k$ with $3\leq k\leq n$ in this case.

{\bf Case 3.} Finally, we discuss the case when $G$ contains at least two vertices $u_i$ and $u_j$ such that\\
(6) $d(u_i)=d(u_j)=\Delta(G)$;\\
(7) both $u_i$ and $u_j$ are in the cycle of $G$;\\
(8) $u_i$ and $u_j$ are not adjacent in $G$.

Then we claim that $px_3(G)\geq \Delta(G)$. Assume to the contrary, $px_3(G)\leq \Delta(G)-1$. Let $c'$ be a $3$-proper coloring of $G$ using colors from $\{1,2,\ldots,\Delta(G)-1\}$. Let the neighbors of $u_i$ be $w_1$, $w_2$, $\ldots$, $w_{\Delta(G)-2}$, $w_{\Delta(G)-1}=u_{i-1}$, $w_{\Delta(G)}=u_{i+1}$, and the neighbors of $u_j$ be $z_1$, $z_2$, $\ldots$, $z_{\Delta(G)-2}$, $z_{\Delta(G)-1}=u_{j-1}$, $z_{\Delta(G)}=u_{j+1}$. Similarly to Case 2, the edges $u_iw_t$ with $t\in \{1,2,\ldots,\Delta(G)-2\}$ are assigned with $\Delta(G)-2$ distinct colors. Without loss of generality, suppose that $c'(u_iw_1)=1$, $c'(u_iw_2)=2$, $\ldots$, $c'(u_iw_{\Delta(G)-2})=\Delta(G)-2$. Thus, either $c'(u_iu_{i-1})=c'(u_iu_{i+1})=\Delta(G)-1$, or there exists at least one edge between $u_iu_{i-1}$ and $u_iu_{i+1}$, say $u_iu_{i-1}$, such that $c'(u_iu_{i-1})=x_1$ with $x_1\in \{1,2,\ldots,\Delta(G)-2\}$. Similarly, the edges $u_jz_{t}$ with $t\in \{1,2,\ldots,\Delta(G)-2\}$ also receive $\Delta(G)-2$ distinct colors. And for the edges $u_ju_{j-1}$ and $u_ju_{j+1}$, either $c'(u_ju_{j-1})=c'(u_ju_{j+1})$, or there exists at least one of them, say $u_ju_{j+1}$, such that $c'(u_ju_{j+1})=c'(u_jz_{x_2})$ with $x_2\in \{1,2,\ldots,\Delta(G)-2\}$.\\
(i) If $c'(u_iu_{i-1})=c'(u_iu_{i+1})$ and $c'(u_ju_{j-1})=c'(u_ju_{j+1})$, then there exists no proper tree containing the vertices $u_{i-1}$, $u_{i+1}$ and $w_1$, a contradiction.\\
(ii) If $c'(u_iu_{i-1})=c'(u_iu_{i+1})$ and $c'(u_ju_{j+1})=c'(u_jz_{x_2})$ with $x_2\in \{1,2,\ldots,\Delta(G)-2\}$, then there exists no proper tree containing the vertices $u_{j+1}$, $u_{j}$ and $z_{x_2}$, a contradiction.\\
(iii) If $c'(u_iu_{i-1})=x_1$ with $x_1\in \{1,2,\ldots,\Delta(G)-2\}$ and $c'(u_ju_{j-1})=c'(u_ju_{j+1})$, then there exists no proper tree containing the vertices $u_{i-1}$, $u_{i}$ and $w_{x_1}$, a contradiction.\\
(iv) If $c'(u_iu_{i-1})=x_1$ with $x_1\in \{1,2,\ldots,\Delta(G)-2\}$ and $c'(u_ju_{j+1})=c'(u_jz_{x_2})$ with $x_2\in \{1,2,\ldots,\Delta(G)-2\}$, then there exists no proper tree containing the vertices $w_{x_1}$, $u_{i-1}$ and $z_{x_2}$, a contradiction.

In summary, we verify that $px_3(G)\geq \Delta(G)$, which deduces that $px_k(G)\geq px_3(G)\geq \Delta(G)$ for $3\leq k\leq n$. Combine with $px_k(G)\leq \Delta(G)$ for $3\leq k\leq n$, we at last arrive at $px_k(G)=\Delta(G)$ for each integer $k$ with $3\leq k\leq n$ in this case.

The proof of this theorem is finished. \qed

We conclude this section with a simple theorem to address the independence of $px_k(G)$ and $rx_k(G)$.

\begin{thm}\label{independence}
For every pair of positive integers $a$, $b$ with $2\leq a\leq b$, there exists a connected graph $G$ such that $px_k(G)=a$ and $rx_k(G)=b$ for each integer $k$ with $3\leq k\leq n$.
\end{thm}
\pf For each pair of positive integers $a$, $b$ with $2\leq a\leq b$, let $G$ be a nontrivial tree with order $n=b+1$ and maximum degree $\Delta(G)=a$. The existence of such a tree is guaranteed by $2\leq a\leq b$. Then based on Theorem \ref{trees} and Lemma \ref{rtrees}, we know that $px_k(G)=\Delta(G)=a$ and $rx_k(G)=n-1=b$ for each integer $k$ with $3\leq k\leq n$. The proof is thus complete. \qed

\section{Graphs with $k$-proper index $n-1$, $n-2$}\label{Characterize}
In this section, we are going to characterize the graphs whose $k$-proper index equals to $n-1$ and $n-2$, respectively, where $3\leq k\leq n$. First of all, we give the following lemma that will be used in the sequel.

\begin{lem}\label{S_n++}
For $n\geq 5$, let $S_n^+$ be the graph obtained by adding a new edge to the $n$-vertices star $S_n$, and $S_n^{++}$ be the graph obtained by adding a new edge to $S_n^+$. Then we have $px_k(S_n^{++})\leq n-3$ for each integer $k$ with $3\leq k\leq n$.
\end{lem}
\pf Let $V(S_n^+)=V(S_n^{++})=\{u,v_1,v_2,\ldots,v_{n-1}\}$. Without loss of generality, set $d_{S_n^+}(u)=d_{S_n^{++}}(u)=n-1$ and $d_{S_n^+}(v_1)=d_{S_n^+}(v_2)=2$. Further, let $e$ be the edge of $S_n^{++}$ added to $S_n^+$. We split the remaining proof into the following two cases depending on the position of $e$.

{\bf Case 1.} The edges $e$ and $v_1v_2$ are vertex-disjoint. Without loss of generality, let $e=v_3v_4$. Then, $G'=G-uv_1-uv_3$ is a spanning tree of $S_n^{++}$ with maximum degree $n-1-2=n-3$. It follows from Theorem \ref{trees} that $px_k(G')=\Delta(G')=n-3$ for $3\leq k\leq n$. Hence, Proposition \ref{subgraph} deduces that $px_k(S_n^{++})\leq px_k(G')=n-3$ for each integer $k$ with $3\leq k\leq n$.

{\bf Case 2.} The edges $e$ and $v_1v_2$ have a common vertex. Without loss of generality, let $e=v_2v_3$. Then, $G''=G-uv_2-uv_3$ is a spanning tree of $S_n^{++}$ with maximum degree $n-1-2=n-3$. Similarly, $px_k(G'')=\Delta(G'')=n-3$ for $3\leq k\leq n$. Hence, we can also get $px_k(S_n^{++})\leq px_k(G'')=n-3$ for each integer $k$ with $3\leq k\leq n$.

Combining the above two cases, now the lemma follows. \qed

\begin{thm}\label{n-1}
Let $G$ be a connected graph of order $n \ (n\geq 4)$. Then for each integer $k$ with $3\leq k\leq n$, $px_k(G)=n-1$ if and only if $G\cong S_n$, where $S_n$ is the star of order $n$.
\end{thm}
\pf Firstly, if $G\cong S_n$, then by Theorem \ref{trees}, we directly obtain $px_k(G)=px_k(S_n)=\Delta(S_n)=n-1$ for $3\leq k\leq n$.

Conversely, suppose $G$ is a graph with $px_k(G)=n-1$ for each integer $k$ with $3\leq k\leq n$. Since $n-1=px_k(G)\leq \Delta(G)$ by Proposition \ref{upperbound2}, meanwhile $\Delta(G)\leq n-1$ holds for any simple graph with order $n$. Then, $\Delta(G)=n-1$. The hypothesis is true if $G\cong S_n$. If $G\ncong S_n$, let $u$ be a vertex of $G$ with $d(u)=\Delta(G)=n-1$. Let $V(G)\backslash u=\{v_1,v_2,\ldots,v_{n-1}\}$ denote the set of the remaining vertices in $G$. Since $G\ncong S_n$, there exist at least two vertices, say $v_1$ and $v_2$, such that they are adjacent in $G$. Set $G'=G[\bigcup\limits_{i=1}^{n-1}uv_i]+v_1v_2$. Then, as $n\geq 4$, $G'$ is a unicyclic graph satisfying the conditions in Case 2 of Theorem \ref{unicyclic} with maximum degree $n-1$. Hence, $px_k(G')=\Delta(G')-1=n-2$ based on the result of Theorem \ref{unicyclic}. Apparently, $G'$ is a spanning subgraph of $G$, therefore $px_k(G)\leq px_k(G')=n-2$ for $3\leq k\leq n$ according to Proposition \ref{subgraph}, contradicting our assumption that $px_k(G)=n-1$. Consequently, $G\cong S_n$.

The proof is thus complete.  \qed

{\bf Remark:} If $G$ is a graph with order $n=3$, then one can check that $px_3(G)=n-1=2$ if and only if $G\cong S_3$ or $G\cong C_3$.

\begin{figure}[h,t,b,p]
\begin{center}
\includegraphics[scale = 0.7]{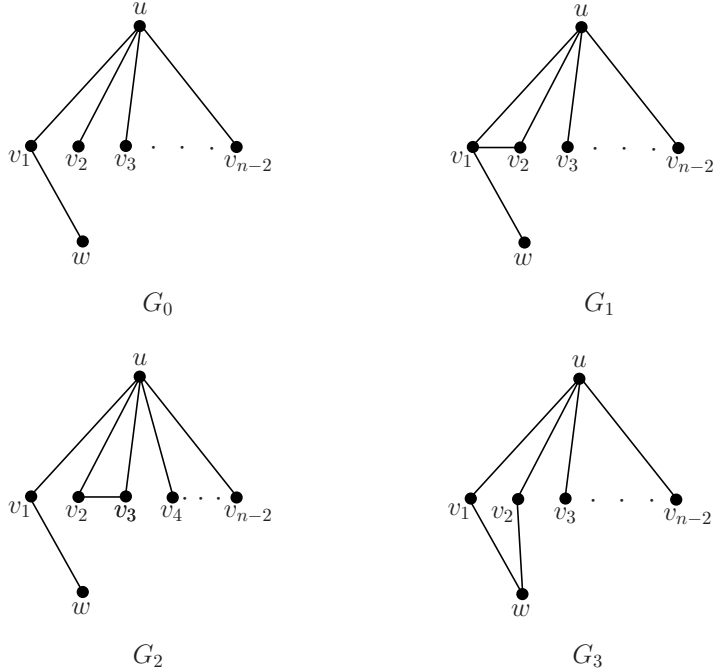}
\caption{The graphs $G_i$ for $i=0,1,2,3$.}\label{fig1}
\end{center}
\end{figure}

\begin{thm}\label{n-2}
Let $G$ be a connected graph of order $n \ (n\geq 5)$. Then for each integer $k$ with $3\leq k\leq n$, $px_k(G)=n-2$ if and only if $G\cong S_n^+$ or $G_0$, where $S_n^+$ is defined in Lemma \ref{S_n++} and $G_0$ is shown in Figure \ref{fig1}.
\end{thm}
\pf On the one hand, if $G\cong S_n^+$, then $G$ is a unicyclic graph with maximum degree $n-1$ satisfying the conditions in Case 2 of Theorem \ref{unicyclic}. Thus $px_k(G)=px_k(S_n^+)=\Delta(G)-1=n-2$ for $3\leq k\leq n$. If $G\cong G_0$, then $G$ is a tree with order $n\geq 5$ and maximum degree $n-2$. Accordingly, by Theorem \ref{trees}, $px_k(G)=px_k(G_0)=\Delta(G)=n-2$ for $3\leq k\leq n$.

On the other hand, if $px_k(G)=n-2$, then by Proposition \ref{upperbound2}, $\Delta(G)\geq px_k(G)=n-2$, which means that $\Delta(G)=n-2$ or $n-1$. The remaining proof is divided into two cases depending on the value of $\Delta(G)$.

{\bf Case 1.} $\Delta(G)=n-1$.\\
In this case, since $px_k(S_n)=n-1$ for $3\leq k\leq n$, as shown before, then $G$ must contain $S_n^+$ as a connected spanning subgraph. If $G\cong S_n^+$, we have known that $px_k(S_n^+)=n-2$ for $3\leq k\leq n$. Now suppose $G\ncong S_n^+$. Then there exists a connected spanning subgraph with the form of $S_n^{++}$ in $G$. Applying Proposition \ref{subgraph} together with Lemma \ref{S_n++}, we arrive at $px_k(G)\leq px_k(S_n^{++})\leq n-3$, a contradiction. Hence, $G\cong S_n^+$ in this case.

{\bf Case 2.} $\Delta(G)=n-2$.\\
Then $G_0$ must be a connected spanning subgraph of $G$. If $G\cong G_0$, then $px_k(G_0)=n-2$ for $3\leq k\leq n$. If $G\ncong G_0$, then there exists at least one edge $e\in E(G)\setminus E(G_0)$. Thus, $G$ contains a connected spanning subgraph isomorphic to $G_1$, $G_2$ or $G_3$, where $G_1$, $G_2$ and $G_3$ are shown in Figure \ref{fig1}. Clearly, one can check that $G_1$, $G_2$ and $G_3$ are all unicyclic graphs with maximum degree $n-2$ satisfying the conditions in Case 2 of Theorem \ref{unicyclic}. Thereupon, by Theorem \ref{unicyclic} as well as Proposition \ref{subgraph}, we directly get that $px_k(G)\leq px_k(G_i)=\Delta(G_i)-1=n-3$ for $3\leq k\leq n$ and $i=1,2$ or $3$, which is a contradiction. Accordingly, $G\cong G_0$ in this case.

In summary, if $px_k(G)=n-2$ for $3\leq k\leq n$, then $G\cong S_n^+$ or $G\cong G_0$. And the proof of this theorem is complete.  \qed

{\bf Remark:} When $n=4$, except for the star $S_4$, other connected graphs with order $4$ are all traceable. Then by Proposition \ref{traceable}, the $k$-proper indices of these graphs equal to $2=n-2$ for each integer $k$ with $3\leq k\leq 4$. While for the star $S_4$, we know that $px_k(S_4)=3$ for $3\leq k\leq 4$. Consequently, we can easily claim that if $G$ is a connected graph of order $n=4$, then for each integer $k$ with $3\leq k\leq 4$, $px_k(G)=n-2=2$ if and only if $G\ncong S_4$.

\end{document}